\documentclass{article}[12pt]
\usepackage{graphicx,pstricks,epsfig,comment}
\usepackage[all]{xy}
\usepackage{subfigure}
\usepackage{rotating}
\usepackage{amsfonts}
\usepackage{amssymb}
\usepackage{amsmath}

\usepackage{indentfirst}
\usepackage{epsfig}
\usepackage{psfrag}
\usepackage{enumerate}
\usepackage{array}
\usepackage{theorem}

\newcommand{\C}{\mathbb{C}}
\newcommand{\R}{\mathbb{R}}

\newcommand{\Z}{\mathbb{Z}}

\newcommand{\Q}{\mathbb{Q}}

\newcommand{\hf}{\frac{1}{2}}
\newcommand{\thf}{\frac{3}{2}}
\newcommand{\sqt}{\frac{\sqrt{3}}{2}}
\newcommand{\phm}{\phantom{-}}
\newcommand{\h}{\hskip0.5em}

\newtheorem{thm}{Theorem}[section]
\newtheorem{cor}[thm]{Corollary}
\newtheorem{lem}{Lemma}[section]
\newtheorem{prop}{Proposition}[section]

{\theorembodyfont{\rmfamily} }
{\theorembodyfont{\rmfamily} }

\newcommand{\Pf}{{\em Proof}. }
\newcommand{\EPf}{\hfill$\Box$\vspace{.5cm}}
\numberwithin{equation}{section}

\title{Complex hyperbolic and projective deformations of small Bianchi groups}

\author{Julien Paupert\footnote{First author partially supported by National Science Foundation Grant DMS-1708463.}, Morwen Thistlethwaite}

\begin{document}
\setlength{\baselineskip}{15pt}

\maketitle

\begin{abstract} The Bianchi groups ${\rm Bi}(d)={\rm PSL}(2,\mathcal{O}_d) < {\rm PSL}(2,\C)$ (where $\mathcal{O}_d$ denotes the ring of integers of $\Q (i\sqrt{d})$, with $d \geqslant 1$ squarefree) can be viewed as subgroups of ${\rm SO}(3,1)$ under the isomorphism ${\rm PSL}(2,\C) \simeq {\rm SO}^0(3,1)$. We study the deformations of these groups into the larger Lie groups ${\rm SU}(3,1)$ and ${\rm SL}(4,\R)$ for small values of $d$. In particular we show that ${\rm Bi}(3)$, which is rigid in ${\rm SO}(3,1)$, admits a 1-dimensional deformation space into ${\rm SU}(3,1)$ and ${\rm SL}(4,\R)$, whereas any deformation of ${\rm Bi}(1)$ into ${\rm SU}(3,1)$ or ${\rm SL}(4,\R)$ is conjugate to one inside ${\rm SO}(3,1)$. We also show that none of the deformations into ${\rm SU}(3,1)$ are both discrete and faithful.
 
\end{abstract}

\section{Introduction}
Let $\Gamma$ be a discrete subgroup of a Lie group $G$, and denote $\iota: \Gamma \rightarrow G$ the inclusion map. 
A \emph{deformation} of $\Gamma$ in $G$ is any continuous 1-parameter family of representations $\rho_t: \Gamma \longrightarrow G$ (for $t$ in some interval $(-\varepsilon,\varepsilon)$) satisfying $\rho_0=\iota$ , and $\rho_t$ not conjugate to $\rho_{0}$ for any $t \neq 0$. We say that $\Gamma$ is \emph{locally rigid} in $G$ if it does not admit any deformations into $G$.

When $G$ is a semisimple real Lie group without compact factors local rigidity of lattices is known to hold in many cases. Weil proved in \cite{W} that $\Gamma$ is locally rigid in $G$ if $G/\Gamma$ is compact and $G$ not locally isomorphic to ${\rm SL}(2,\R)$. Garland and Raghunathan extended this result to the case where $\Gamma$ is a non-cocompact lattice in a rank-1 semisimple group $G$ not locally isomorphic to ${\rm SL}(2,\R)$ or ${\rm SL}(2,\C)$ (Theorem 7.2 of \cite{GR}). 

Mostow's strong rigidity theorem (\cite{Mo}) is both stronger and less general. It asserts that the inclusion map of any lattice $\Gamma$ in a rank-1 semisimple Lie group $G$, not locally isomorphic to ${\rm SL}(2,\R)$, is the unique faithful representation of $\Gamma$ into $G$ whose image is a lattice in $G$ (up to conjugation). In particular, $\Gamma$ has no deformations which are discrete and faithful with image a lattice in $G$.  

Both local rigidity and Mostow rigidity fail dramatically in ${\rm SL}(2,\R)$. 
Most lattices in $G={\rm SL}(2,\R)$ are well known to have large deformation spaces (of lattice embeddings) inside ${\rm SL}(2,\R)$, which are known as the Teichm\"uller space of the corresponding lattice or quotient orbifold. For example the fundamental group of a compact surface of genus $g \geqslant 2$ enjoys a real $(6g-6)$-dimensional Teichm\"uller space. In the context of this paper it is also worth noting that some ``small'' lattices in ${\rm SL}(2,\R)$ such as compact triangle groups are in fact locally rigid in ${\rm SL}(2,\R)$. 

In the case of $G={\rm SL}(2,\C)$, all lattices satisfy Mostow rigidity and cocompact lattices satisfy local rigidity by Weil's above result. Thurston showed however that many cusped lattices in ${\rm SL}(2,\C)$ admit rich deformation spaces into ${\rm SL}(2,\C)$, the so-called \emph{Dehn surgery spaces}. More specifically (see Theorem 5.8.2 of \cite{T}), if $M$ is a finite-volume complete orientable hyperbolic 3-manifold with $k \geqslant 1$ cusps, then $\pi_1(M)$ admits a (real) $2k$-dimensional deformation space into ${\rm SL}(2,\C)$. (Here the original lattice embedding $\pi_1(M) \rightarrow {\rm SL}(2,\C)$ is the holonomy representation of the unique complete, finite-volume hyperbolic structure on $M$). These deformations are not both discrete and faithful, but there exist infinitely many deformations with discrete and finite-covolume image in any neighborhood of the lattice embedding (corresponding to manifolds obtained by \emph{Dehn surgery} on the cusps of $M$).   

Dunbar and Meyerhoff showed in \cite{DM} that Thurston's results carry through to the case of (complete, finite-volume, orientable) cusped hyperbolic 3-orbifolds, but only for those cusps whose cross-sections (which are in general Euclidean 2-orbifolds) are a torus or a $(2,2,2,2)$-pillow, i.e. a sphere with four cone angles of angle $\pi/2$. More specifically, let $Q$ be such an orbifold, denote $\Gamma=\pi_1^{orb}(Q)$ a lattice in ${\rm SL}(2,\C)$ such that $Q={\rm H}^3_\R/\Gamma$ and let $k$ be the number of cusps of $Q$ whose cross-sections are either tori or $(2,2,2,2)$-pillows. Theorem~5.3 of \cite{DM} implies that the deformation space of $\Gamma$ into ${\rm SL(2,\C)}$ contains a smooth (real) $2k$-dimensional family of deformations. We will abbreviate this as: ${\rm dim}_{[\rho_{{\rm hyp}}]} \, \chi (\Gamma,{\rm SL}(2,\C)) \geqslant 
2k$, where $\chi (\Gamma,G)$ is the \emph{character variety} of $\Gamma$ into $G={\rm SL}(2,\C)$ and $\rho_{{\rm hyp}}: \Gamma \rightarrow G$ is the lattice embedding (unique up to conjugation by Mostow), even if $\chi (\Gamma,G)$ is not necessarily smooth at $\rho_{{\rm hyp}}$. (See Section~\ref{varieties} for a precise definition of $\chi (\Gamma,G)$ and more details on character varieties and smoothness).

In this paper we focus on the \emph{Bianchi groups} ${\rm Bi}(d)={\rm PSL}(2,\mathcal{O}_d) < {\rm PSL}(2,\C)$, where $\mathcal{O}_d$ denotes the ring of integers of $\Q (i\sqrt{d})$, with $d \geqslant 1$ squarefree. These are well known to be (arithmetic) lattices in  ${\rm PSL}(2,\C)$, with number of cusps equal to the class number $h_d$ of  $\Q (i\sqrt{d})$. In particular, ${\rm Bi}(d)$ has a single cusp when $\mathcal{O}_d$ is Euclidean, namely when $d=1,2,3,7,11$. (${\rm Bi}(d)$ also has a single cusp when $d=19,43,67,163$ but we will not consider those here). It is also known that the cusp cross-sections of ${\rm Bi}(d)$ are all tori when $d \neq 1,3$, and are a $(2,2,2,2)$-pillow when $d=1$ and a $(3,3,3)$-turnover (sphere with three cone points of angle $\pi/3$) when $d=3$. See \cite{MR} and \cite{F} for more details on Bianchi groups and section 2.4 of \cite{BR} for a discussion of the Dunbar-Meyerhoff result in the case of the Bianchi groups. It follows therefore that:

\begin{thm} With the above notation: 
\begin{itemize}
\item ({\bf [DM]}) ${\rm dim}_{[\rho_{{\rm hyp}}]} \, \chi ({\rm Bi}(d),{\rm PSL}(2,\C)) \geqslant 2h_d$ for all $d \neq 3$.
\item ({\bf [Se]}) ${\rm dim}_{[\rho_{{\rm hyp}}]} \, \chi ({\rm Bi}(3),{\rm PSL}(2,\C)) =0$.
\end{itemize} 
\end{thm}

The second item follows from the fact that ${\rm Bi}(3)$ is locally rigid in ${\rm PSL}(2,\C)$, since following Serre (\cite{Se}) it has property (FA). In fact, more recently it was shown in \cite{BMRS} that ${\rm Bi}(3)$ has the much stronger \emph{Galois rigidity} property, namely that the only two Zariski-dense repesentations of ${\rm Bi}(3)$ into ${\rm PSL}(2,\C)$ are (up to conjugation) the lattice embedding and its complex conjugate.

Under the isomorphism ${\rm PSL}(2,\C) \simeq {\rm SO}^0(3,1)$, one may also consider the Bianchi groups as subgroups of ${\rm SO}(3,1)$ and ask whether the lattice embedding $\rho_{{\rm hyp}}: {\rm Bi}(d) \longrightarrow {\rm SO}(3,1)$ admits any further deformations into other Lie groups $G$ containing ${\rm SO}(3,1)$. In this paper we consider \emph{complex hyperbolic} and \emph{projective} deformations of the Bianchi groups, that is deformations into the Lie groups $G={\rm SU}(3,1)$ and ${\rm SL}(4,\R)$ respectively. Denote $\iota$ the inclusion ${\rm SO}(3,1) \longrightarrow {\rm SU}(3,1)$. By analogy with the 2-dimensional \emph{$\R$-Fuchsian} terminology, we will call $\iota \circ \rho_{\rm hyp}$ the \emph{$\R$-Kleinian embedding} of $\Gamma$ into ${\rm SU}(3,1)$.

Projective deformations and complex hyperbolic deformations of real hyperbolic lattices/discrete groups are known to be equivalent in the following sense (Theorem~2.2 of \cite{CLT1}):

\begin{thm}[\cite{CLT1}] Let $\Gamma$ be a finitely generated group, and let $\rho:\Gamma \longrightarrow {\rm SO}^0(n,1)$ be a smooth point of the representation variety ${\rm Hom}(\Gamma, {\rm SL}(n+1,\R))$. Then $\iota \circ \rho$ is also a smooth point of ${\rm Hom}(\Gamma, {\rm SU}(n,1))$, and near $\rho$ the real dimensions of ${\rm Hom}(\Gamma, {\rm SL}(n+1,\R))$ and ${\rm Hom}(\Gamma, {\rm SU}(n,1))$ are equal.
\end{thm}

Our main results are the following; in the next statement one can replace ${\rm SU}(3,1)$ with ${\rm SL}(4,\R)$ in view of the above result. 

\begin{thm}\label{maindim} For $d=3,1$ the character variety $\chi ({\rm Bi}(d),{\rm SU}(3,1))$ (resp. $\chi ({\rm Bi}(d),{\rm SL}(4,\C))$) is smooth at $[\iota \circ \rho_{{\rm hyp}}]$, with:
\begin{itemize}
\item ${\rm dim}_{[\iota \circ \rho_{{\rm hyp}}]} \, \chi ({\rm Bi}(3),{\rm SU}(3,1)) =1$ (resp. 2) 
\item ${\rm dim}_{[\iota \circ \rho_{{\rm hyp}}]} \, \chi ({\rm Bi}(1),{\rm SU}(3,1)) =2$ (resp. 4) 
\end{itemize} 
\end{thm}

We also compute the dimension of the Zariski tangent spaces $H^1({\rm Bi}(d),{\mathfrak sl}(4,\R))_{[\iota \circ \rho_{{\rm hyp}}]}$ when $d=2,5,6,7,11,$ $15,19$ in Proposition~\ref{Zariski}, but we do not know whether $[\iota \circ \rho_{{\rm hyp}}]$ is a smooth point in those cases.

Porti studied in \cite{Po} deformations of general hyperbolic 3-orbifold groups $\Gamma$ into complex semi-simple Lie groups $G$ including ${\rm SL}(n,\C)$, where the original representation is obtained by post-composing the holonomy representation with the irreducible representation ${\rm SL}(2,\C) \rightarrow {\rm SL}(n,\C)$ (so, analogous to the \emph{Fuchsian} representation defining the \emph{Hitchin} component when $G={\rm SL}(n,\R)$, see e.g. \cite{ALS} and \cite{LT2} for hyperbolic 2- and 3-orbifolds). When $\Gamma$ is a Bianchi group and $G={\rm SL}(4,\C)$, Porti's results apply to a point in a different component of the character variety than the one we consider here. In particular, he shows that $\chi ({\rm Bi}(d),{\rm SL}(4,\C))$ is smooth at that point, with dimension three times the number of cusps when all boundary components are tori, that is $3h_d$ when $d \neq 1,3$. One can show that the $\R$-Kleinian embedding $\iota \circ \rho_{{\rm hyp}}$ we consider is \emph{strongly regular} on the boundary tori, so that by Theorem 1.2 of \cite{Po}, $\chi ({\rm Bi}(d),{\rm SL}(4,\C))$ has dimension at least $3h_d$ at the $\R$-Kleinian embedding (see also \cite{FG} in the manifold case). However we suspect that the dimension is strictly larger than this for higher values of $d$, in particular the 
dimension of the Zariski tangent spaces that we compute is strictly greater when $d \geqslant 7$. It is interesting to note that for certain manifold subgroups, namely the figure-eight knot group in ${\rm Bi}(3)$ (respectively the Whitehead link group in ${\rm Bi}(1)$), the $\R$-Kleinian embedding is in fact a smooth point with dimension 3 (resp. 6) by results of Ballas--Danciger--Lee, see Section~\ref{manifolds}.

From the point of view of geometric structures, discrete and faithful representations are of particular interest. Indeed they correspond to \emph{uniformizable} (or \emph{complete}) geometric structures, that is those structures whose associated developing map is a covering map, in the case where the model space carries an invariant Riemannian metric. More specifically, given a pair $(G,X)$ with $G$ a Lie group acting transitively by isometries on a Riemannian manifold $X$, for any manifold $M$ there is a bijection between the space of (marked) complete $(G,X)$-structures on $M$ and the \emph{discrete and faithful character variety} $\chi_{DF}(\pi_1(M),G)$, which is the space of conjugacy classes of discrete and faithful representations of $\pi_1(M)$ into $G$. (See e.g. section 8.1 of \cite{Gol2} for more details). 

From this point of view the lattice embedding $\rho_{{\rm hyp}}: \Gamma \rightarrow {\rm SO}(3,1)$ is the holonomy representation corresponding to the (unique) complete, finite-volume hyperbolic structure on the orbifold  $Q={\rm H}^3_\R/\Gamma$. Moreover it is well known (among experts) that the $\R$-Kleinian embedding $\iota \circ \rho_{\rm hyp} : \Gamma \rightarrow {\rm SU}(3,1)$ is the holonomy of a complete, complex hyperbolic structure on the tangent (orbi-)bundle of $Q$, in other words that ${\rm H}_\C^3/\iota(\rho_{\rm hyp}(\Gamma))$ is diffeomorphic to the tangent (orbi-)bundle of $Q$; see section~\ref{structures} for more details.

\begin{thm}\label{mainDF} For $d=1,3$, $[\iota \circ \rho_{{\rm hyp}}]$ is the only point in its irreducible component of $\chi ({\rm Bi}(d),{\rm SU}(3,1))$ comprising discrete and faithful representations. In particular, the complete complex hyperbolic structure on the tangent orbibundle to the corresponding Bianchi orbifold is locally rigid.
\end{thm}

It is well known that the modular surface ${\rm H}^2_\R/{\rm PSL}(2,\Z)$ is embedded in the Bianchi orbifold ${\rm H}^3_\R/{\rm Bi}(d)$ for all $d\neq1,3$ (see e.g. section 6.3.2 of \cite{F} ), so that ${\rm Bi}(d)$ admits a 1-parameter family of \emph{bending deformations} into ${\rm SL}(4,\R)$, in the sense of Johnson-Millson (\cite{JM}).
Therefore there is at least one dimension worth of smooth deformations of ${\rm Bi}(d)$ into ${\rm SL}(4,\R)$, transverse to the Dehn surgery deformations of Dunbar--Meyerhoff, for all $d \neq 1,3$. In particular, when $d=7$ this gives (assuming smoothness at $[\rho_{{\rm hyp}}]$ for notational convenience): ${\rm dim}_{[\rho_{{\rm hyp}}]} \, \chi ({\rm Bi}(7),{\rm SL}(4,\R)) \geqslant 3$.

The Bianchi groups ${\rm Bi}(1)$ and ${\rm Bi}(3)$ are also well known to be commensurable to several Coxeter groups (see e.g. Figures 13.2 and 13.3 of \cite{MR}); in particular ${\rm Bi}(1)$ (resp. ${\rm Bi}(3)$) is the index-2 orientation-preserving subgroup of a Coxeter group $\Gamma_1$ (resp. $\Gamma_3$), whose diagrams are pictured in Figures 1 and 2. Projective deformations of tetrahedral Coxeter groups were studied in \cite{Ma}, and it follows from the main result there that the corresponding Coxeter groups $\Gamma_1$ (resp. $\Gamma_3$) have no deformations (resp. a smooth 1-parameter family of deformations) into ${\rm SL}(4,\R)$. (According to Section 7.4 of \cite{BDL}, the result for $\Gamma_3$ also follows from \cite{Ben}). This shows in particular that ${\rm dim}_{[\rho_{{\rm hyp}}]} \, \chi ({\rm Bi}(3),{\rm SU}(3,1)) \geqslant 1$; our result shows that any deformation of ${\rm Bi}(3)$ into ${\rm SU}(3,1)$ or ${\rm SL}(4,\R)$ is induced by a deformation of the Coxeter group $\Gamma_3$ (at least, in the irreducible component of $[\rho_{{\rm hyp}}]$).

The paper is organized as follows. In Section 2 we review generalities about representation and character varieties and in Section 3 Swan's presentations of small Bianchi groups. In Section 4 we specialize to deformations of Bianchi groups into ${\rm SU}(3,1)$ and ${\rm SL}(4,\R)$, and compute the corresponding Zariski-tangent spaces for $d=2,5,6,7,11,15,19$ and explicit deformations of ${\rm Bi}(3)$, proving the results announced above. In Section 5 we review the geometric structures associated to $\R$-Kleinian embeddings in ${\rm SU}(3,1)$, and results about deformations of manifold subgroups such as the figure-eight knot and Whitehead link groups.

The first author would like to thank Martin Deraux, Pierre Will, and Alan Reid for helpful discussions about some of the material discussed here.

\begin{figure}\label{fig1}
  \centering
  \caption{The Coxeter group $\Gamma_1$ containing ${\rm Bi}(1)$ with index 2}
  \epsfig{figure=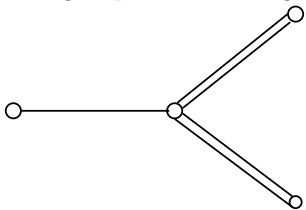, width=4cm}
\end{figure}

\begin{figure}\label{fig2}
  \centering
  \caption{The Coxeter group $\Gamma_3$ containing ${\rm Bi}(3)$ with index 2}
  \epsfig{figure=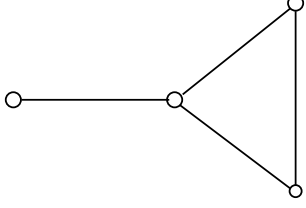, width=4cm}
\end{figure}

\section{Representation varieties, character varieties and smoothness}\label{varieties}

Let $\Gamma$ be a finitely generated group and $G$ a linear algebraic group, more specifically a $\Q$-algebraic subgroup of ${\rm GL}(n,\R)$, or the set of real points in ${\rm GL}(n,\C)$ of a $\Q$-algebraic group (e.g. $G={\rm SL}(n,\R), {\rm SO}(n,1)$ or ${\rm SU}(n,1)$).

The \emph{representation variety} ${\rm Hom}(\Gamma,G)$ is the space of all representations of $\Gamma$ into $G$.
Since $G$ is a linear algebraic group and $\Gamma$ is finitely generated, ${\rm Hom}(\Gamma,G)$ is an affine algebraic variety. 

The \emph{character variety} $\chi (\Gamma,G)$ is the (maximal Hausdorff quotient of) the topological quotient ${\rm Hom}(\Gamma,G)/G$ of the representation variety under $G$ (acting on representations by postcomposition by inner automorphisms). The latter is in fact also an algebraic variety, and can be alternatively described as the quotient in the sense of geometric invariant theory; see \cite{A} for more details on character varieties in real forms of ${\rm SL}(n,\C)$.

\subsection{Infinitesimal deformations and cohomology}

We start with an elementary computation to motivate the definition of the relevant cohomology groups due to Weil (\cite{W}); see \cite{Bes} or section 3.1 of \cite{HP} for more details. Given a finitely generated group $\Gamma$ and a linear Lie group $G$ as above, assume that $\rho_t:\Gamma \longrightarrow G$ (with $t \in (-\varepsilon, \varepsilon)$ for some $\varepsilon>0$) is a $C^1$-family of representations, in the sense that for all $\gamma \in \Gamma$, the map $t \mapsto \rho_t(\gamma)$ is $C^1$. Then the tangent vectors to the various $\rho_t (\gamma)$ ($\gamma \in \Gamma$), appropriately (right-) translated into the Lie algebra $\mathfrak{g}$ will satisfy the following cocycle condition.

{\bf Cocycles:} Denote, for any $\gamma \in \Gamma$ and $t \in (-\varepsilon,\varepsilon)$, $d_{\rho_0}(t,\gamma)=\rho_t(\gamma)\rho_0(\gamma)^{-1}$ and $h_{\rho_0}(\gamma)=\frac{d}{dt}\bigr\rvert_{t=0}d_{\rho_0}(t,\gamma) \in \mathfrak{g}$. Then, for any $\gamma,\gamma' \in \Gamma$ we have:

$$\begin{array}{rcl}
d_{\rho_0}(t,\gamma\gamma') & = & \rho_t(\gamma\gamma')\rho_0(\gamma\gamma')^{-1} \\
 & = & \rho_t(\gamma)\rho_t(\gamma')\rho_0(\gamma')^{-1}\rho_0(\gamma)^{-1} \\
 & = & \left[\rho_t(\gamma)\rho_0(\gamma)^{-1}\right] \left[ \rho_0(\gamma)\rho_t(\gamma')\rho_0(\gamma')^{-1}\rho_0(\gamma)^{-1}\right] \\
 & = & d_{\rho_0}(t,\gamma)\,  \rho_0(\gamma) d_{\rho_0}(t,\gamma')\rho_0(\gamma)^{-1} 
\end{array}
$$
Therefore, taking derivatives and evaluating at $t=0$ we get:
$$\begin{array}{rcl}
h_{\rho_0}(\gamma\gamma') & = & h_{\rho_0}(\gamma)+\rho_0(\gamma) h_{\rho_0}(\gamma')\rho_0(\gamma)^{-1} \\
& = &  h_{\rho_0}(\gamma)+{\rm Ad}\, \rho_0(\gamma)\cdot h_{\rho_0}(\gamma'),
\end{array}
$$
where:
$$\begin{array}{rccl}
{\rm Ad}\, \rho_0 : & \Gamma & \longrightarrow & {\rm GL}(\mathfrak{g}) \\
 & \gamma & \longmapsto & [ v \mapsto \rho_0(\gamma)v\rho_0(\gamma)^{-1} ]
 \end{array}
 $$

{\bf Deformations by conjugation:} As a special case, if $(g_t)$ ($t \in (-\varepsilon,\varepsilon)$) is a $C^1$-path in $G$ with $g_0=e$, we obtain a $C^1$-deformation $(c_t)$ of $\rho_0$ by defining: $c_t(\gamma)=g_t \rho_0(\gamma)g_t^{-1}$ for any $\gamma \in \Gamma$. Then:
$$ d_{c_0}(t,\gamma)=g_t\rho_0(\gamma)g_t^{-1}\rho_0(\gamma)^{-1}, \ {\rm and}$$
$$h_{c_0}(\gamma) = X-\rho_0(\gamma)X\rho_0(\gamma)^{-1}=X-{\rm Ad}\, \rho_0(\gamma)\cdot X, $$
where: $X = \frac{d}{dt}\bigr\rvert_{t=0}g_t \in \mathfrak{g}$. This motivates the following definitions:

$$ {\rm (Cocycles)} \ \ Z^1(\Gamma,\mathfrak{g})_\rho = \left\lbrace u_1: \Gamma \longrightarrow \mathfrak{g} \, | \, (\forall \gamma,\delta \in \Gamma) \, u_1(\gamma\delta)=u_1(\gamma) +  {\rm Ad}\, \rho(\gamma)\cdot u_1(\delta)    \right\rbrace,$$
$$ {\rm (Coboundaries)} \ \ B^1(\Gamma,\mathfrak{g})_\rho = \left\lbrace v_1: \Gamma \longrightarrow \mathfrak{g} \, | \, (\exists X \in \mathfrak{g})(\forall \gamma \in \Gamma) \, v_1(\gamma)=X -  {\rm Ad}\, \rho(\gamma)\cdot X    \right\rbrace,$$
$$ {\rm (Cohomology \ groups)} \ \  H^1(\Gamma,\mathfrak{g})_\rho= Z^1(\Gamma,\mathfrak{g})_\rho/ B^1(\Gamma,\mathfrak{g})_\rho.$$

\subsection{Computing Zariski tangent spaces}\label{zariski}

See Section 2 of \cite{CLT1} for a more detailed discussion. Assuming for simplicity that $\Gamma$ is finitely presented and $G={\rm SL}(n,\R)$, the affine variety ${\rm Hom}(\Gamma,G)$ has the following concrete description. Consider a finite presentation $\langle g_1,...,g_m \, | \, w_1(g_1,...,g_m)=...=w_k(g_1,...,g_m)=1 \rangle$ of $\Gamma$, with $w_1,...,w_k$ words in the letters $g_1,...,g_m$ and their inverses. Then a representation $\rho: \Gamma \longrightarrow G$ is determined by an $m$-tuple $(A_1,...,A_m) \in G^m$ satisfying the equations $w_1(A_1,...,A_m)=...=w_k(A_1,...,A_m)={\rm Id}$. Moreover, since $G={\rm SL}(n,\R) \subset {\rm M}(n,\R)$, identifying ${\rm M}(n,\R)$ with $\R^{n^2}$ will identify ${\rm Hom}(\Gamma,G)$ with $V=f^{-1}(\{0\})$, where:

$$\begin{array}{rrcl}
f: & (\R^{n^2})^m & \longrightarrow & \R^m \times (\R^{n^2})^k \\
 & (A_1,...,A_m) & \longmapsto & ({\rm det}\, A_1-1, ..., ,{\rm det}\, A_m-1, w_1(A_1,...,A_m)-{\rm Id},..., w_k(A_1,...,A_m)-{\rm Id})
 \end{array}
$$
From differential geometry we know that, if $f$ is a submersion, $V$ is smooth with tangent space $T_pV={\rm Ker} \, df_p$ for any $p \in V$. However this is often not the case, in particular it can never happen if $k \geqslant m$, which will be the case for all the presentations we use in this paper.  (A word about "expected dimension?")

The linear subspace $T^Z_pV={\rm Ker} \, df_p \subset (\R^{n^2})^m$ is usually called the \emph{Zariski tangent space} to $V$ at $p$; by construction it contains the tangent vector at $p$ to any smooth curve through $p$ in $V$, hence its dimension is an upper bound for the dimension of any smooth irreducible component of $V$ containing $p$. 

It is well known that $T^Z_\rho{\rm Hom}(\Gamma,G) \simeq Z^1(\Gamma,\mathfrak{g})_\rho$ as real vector spaces (see e.g. p. 152 of \cite{W}). In particular, if the conjugation orbit $C_{\rho}=\{ g \rho g^{-1} \, | \, g \in G \} \subset {\rm Hom}(\Gamma,G)$ is smooth at $\rho$ then: $H^1(\Gamma,\mathfrak{g})_\rho \simeq T^Z_\rho{\rm Hom}(\Gamma,G) / TC_\rho$.

\subsection{Smoothness}
See Section 3.1 of \cite{BDL} for a more detailed discussion. 
\begin{lem}\label{dimH1} Let $\Gamma$ be a finitely generated group, $G$ a Lie group and $\rho: \Gamma \longrightarrow G$  an irreducible representation. Then the conjugation orbit $C_{\rho}=\{ g \rho g^{-1} \, | \, g \in G \} \subset {\rm Hom}(\Gamma,G)$ is smooth at $\rho$ with dimension ${\rm dim}\, G$. In particular, ${\rm dim}\, H^1(\Gamma,\mathfrak{g})_\rho ={\rm dim}\, T^Z_\rho{\rm Hom}(\Gamma,G) -{\rm dim}\, G$. 
\end{lem}
We will give concrete examples of these computations in Section~\ref{zariski2}.
The concrete criterion we will use to determine smoothness of the character variety is the following:

\begin{prop} Let $\Gamma$ be a finitely generated group, $G$ a Lie group and $\rho: \Gamma \longrightarrow G$  an irreducible representation. Denote $k={\rm dim}\, H^1(\Gamma,\mathfrak{g})_\rho$ and assume that there exist a neighborhood $U$ of $0$ in $\R^k$ and a smooth family of representations $\rho_t: \Gamma \longrightarrow G$ $(t \in U)$ such that $\rho_0=\rho$ and $\rho_t$ not conjugate to $\rho_{t'}$ for any $t \neq t' \in U$. Then $\chi(\Gamma,G)$ is smooth at $\rho$ with dimension $k$.
\end{prop}

\section{Bianchi groups}

For any squarefree $d \geqslant 1$ we consider the Bianchi group ${\rm Bi}(d)={\rm PSL}(2,\mathcal{O}_d) < {\rm PSL}(2,\C)$, where $\mathcal{O}_d$ denotes as usual the ring of integers of $\Q (i\sqrt{d})$. Recall that $\mathcal{O}_d=\Z[\tau]$, where $\tau=i\sqrt{d}$ if $d=1,2$ mod 4 and $\tau=\frac{1+i\sqrt{d}}{2}$ if $d=3$ mod 4. Swan found explicit presentations for most Bianchi groups ${\rm Bi}(d)$ with $d \leqslant 19$ in \cite{Sw}. We now list Swan's generators and presentations, for the 1-cusped groups with $d=1,2,3,7,11,19$ and the 2-cusped groups with $d=5,6,15$  (when $d=3$ we use $\tau =\frac{-1+i\sqrt{3}}{2}$ for these generators). The generators include for all $d$:

$$\begin{array}{ccc}
T=\left(\begin{array}{cc}
1 & 1 \\
0 & 1 \end{array}\right)
&
U=\left(\begin{array}{cc}
1 & \tau \\
0 & 1 \end{array}\right)
&
A=\left(\begin{array}{cc}
0 & -1 \\
1 & 0 \end{array}\right)
\end{array}
$$
as well as, when $d=1,3$:
$$L=\left(\begin{array}{cc}
\tau^{-1} & 0 \\
0 & \tau \end{array}\right),
$$
when $d=19$:
$$B=\left(\begin{array}{cc}
1-\tau & 2 \\
2 & \tau \end{array}\right),
$$
when $d=15$:
$$C=\left(\begin{array}{cc}
4 & 1-2\tau  \\
2\tau -1 & 4 \end{array}\right),
$$
when $d=5$:
$$ \begin{array}{ccc}
B=\left(\begin{array}{cc}
-\tau & 2  \\
2 & \tau \end{array}\right)
&  {\rm and} & C=\left(\begin{array}{cc}
-\tau -4 & 2\tau  \\
2\tau  & \tau -4 \end{array}\right),
\end{array}
$$
and when $d=6$:
$$ \begin{array}{ccc}
B=\left(\begin{array}{cc}
-1-\tau & 2-\tau  \\
2 & 1+ \tau \end{array}\right)
&  {\rm and} & C=\left(\begin{array}{cc}
5 & -2\tau  \\
2\tau  & 5 \end{array}\right).
\end{array}
$$

The presentations given by Swan are the following:

$$
\begin{array}{c}
{\rm Bi}(1)= \langle T,U,L,A \, | \, [T,U]=L^2=(TL)^2=(UL)^2=(AL)^2=A^2=(TA)^3=(UAL)^3=1 \rangle \\
\\
{\rm Bi}(3)= \langle T,U,L,A \, | \, [T,U]=L^3=(AL)^2=A^2=(TA)^3=(UAL)^3=1, \, L^{-1}UL=T, \, L^{-1}TL=T^{-1}U^{-1} \rangle \\ 
\\
{\rm Bi}(2)= \langle T,U,A \, | \, [T,U]=A^2=(TA)^3=(AU^{-1}AU)^2=1 \rangle \\
\\
{\rm Bi}(7)= \langle T,U,A \, | \, [T,U]=A^2=(TA)^3=(ATU^{-1}AU)^2=1 \rangle \\
\\
{\rm Bi}(11)= \langle T,U,A \, | \, [T,U]=A^2=(TA)^3=(ATU^{-1}AU)^3=1 \rangle \\
\\
{\rm Bi}(19)= \langle T,U,A,B \, | \, [T,U]=A^2=(TA)^3=B^3=(BT^{-1})^3=(AB)^2=(AT{-1}UBU^{-1})^2=1 \rangle \\
\\
{\rm Bi}(15)= \langle T,U,A,C \, | \, [T,U]=[A,C]=A^2=(TA)^3=1, UCUAT=TAUCU \rangle \\
\\
{\rm Bi}(5)= \langle T,U,A,B,C \, | \, [T,U]=A^2=(TA)^3=B^2=(AB)^2=(ATUBU^{-1})^2=1, ACA=TCT^{-1}=UBU^{-1}CB \rangle \\
\\
{\rm Bi}(6)= \langle T,U,A,B,C \, | \, [T,U]=[A,C]=A^2=(TA)^3=B^2=(ATB)^3=(ATUBU^{-1})^3=1, CTUB=TBCU \rangle 
\end{array}
$$

\section{Deformations of Bianchi groups into ${\rm SU}(3,1)$}

\subsection{Infinitesimal deformations}\label{zariski2}

We first compute the spaces of infinitesimal deformations $H^1({\rm Bi}(d),{\mathfrak sl}(4,\R))_{\rho_0}$.
The principle of the computation is very simple, but requires computing the rank of rather large matrices (from 48 by 67 to 80 by 133 in the cases considered here) so are best handled with formal computation software such as Mathematica or Maple. We write the details of the computation when $d=7$; the others are similar but may have more generators and/or relations.

We will relabel the generators in Swan's presentation for ${\rm Bi}(7)$ above as $a,b,c$ for convenience, and denote their inverses by $A,B,C$. Swan's presentation now reads:
$$
{\rm Bi}(7)= \langle a,b,c \, | \, [a,b]=c^2=(ac)^3=(caBcb)^2=1 \rangle
$$
To decrease degrees, since we are considering linear representations, we may rewrite the relations as: ${\rm Rel}1:\, ab-ba$, ${\rm Rel}2:\, c^2 -{\rm Id}$, ${\rm Rel}3:\, aca-CAC$, ${\rm Rel}4: \,caBcb-BCbAC$.
We start with the following generators in ${\rm SO}(3,1)$:

$$
\begin{array}{cc}
a_0=\left(\begin{array}{cccc}
3/2 & -1/2 & 1 & 0 \\
1/2 & 1/2 & 1 & 0 \\
1 & -1 & 1 & 0 \\
0 & 0 & 0 & 1 
\end{array}\right),
&
b_0=\left(\begin{array}{cccc}
2 & -1 & 1/2 & \sqrt{7}/2 \\
1 & 0  & 1/2 & \sqrt{7}/2 \\
1/2 & -1/2 & 1 & 0 \\
 \sqrt{7}/2 & - \sqrt{7}/2 & 0 & 1 
\end{array}\right),
\end{array}
$$

$$
c_0=\left(\begin{array}{cccc}
1 & 0 & 0 & 0 \\
0 & -1 & 0 & 0 \\
0 & 0 & -1 & 0 \\
0 & 0 & 0 & 1 
\end{array}\right).
$$

Consider the polynomial map:
$$\begin{array}{rrcl}
{\rm Eqs}: & (\R^{16})^3 & \longrightarrow & \R^3 \times (\R^{16})^4 \\
 & (a,b,c) & \longmapsto & ({\rm det}\, a-1, {\rm det}\, b-1, {\rm det}\, c-1, {\rm Rel}1,...,{\rm Rel}4)
 \end{array}
$$
(Taking the inverse of a matrix is only a rational map, but since we are requiring the determinants to be 1 it is in fact polynomial in this case). 

The above presentation then identifies ${\rm Hom}({\rm Bi}(7),{\rm SL}(4,\R))$ with ${\rm Eqs}^{-1}(\{(0,...,0)\})$, and the Zariski tangent space to ${\rm Hom}({\rm Bi}(7),{\rm SL}(4,\R)$ at the point $(a_0,b_0,c_0)$ is identified with the kernel of the differential $d_{(a_0,b_0,c_0)} {\rm Eqs}$. Coding the map ${\rm Eqs}$ and computing the rank of $d_{(a_0,b_0,c_0)} {\rm Eqs}$ (a real 48 by 67 matrix) is routine in any formal software system; Maple and Mathematica tell us that in the present case the rank is 30, hence the dimension of the kernel is $48-30=18$. Now we conclude by Lemma~\ref{dimH1} that the dimension of $H^1 ({\rm Bi}(7),{\rm SL}(4,\R))_{\rho_0}$ is $18-{\rm dim}\, {\rm SL}(4,\R)=18-15=3$. Similar computations give the following: 
\begin{prop}\label{Zariski} The space of infinitesimal deformations $H^1({\rm Bi}(d),{\mathfrak sl}(4,\R))_{\rho_0}$ has dimension:
\begin{itemize}
\item 1 when $d=3$,
\item 2 when $d=1$,
\item 3 when $d=2$ or $7$,
\item 4 when $d=11$,
\item 5 when $d=19$,
\item 7 when $d=5$ or $15$,
\item 9 when $d=6$.
\end{itemize} 
\end{prop}

\subsection{Computing the deformation spaces}

In \cite{CLT2} a method was introduced for exact computation of character
varieties of fundamental groups of low-dimensional manifolds and orbifolds.
This method is described in \cite{CLT2} and \cite{LT1} in some detail; it
often works well when the variety has dimension at most $2$, this being
the case for many interesting examples.  For this project we computed the
full deformation spaces for ${\rm Bi}(1)$ and ${\rm Bi}(3)$; we also computed
a $2$--dimensional subvariety of the $3$--dimensional deformation space
for ${\rm Bi}(7)$.

Here is a description of how this technique applies to the $1$--dimensional
deformation space for ${\rm Bi}(3)$.

The starting point is the holonomy representation $\rho_{\rm hyp}$ of ${\rm Bi}(3)$ into the
Lorentz group ${\rm O}(3,1)$, namely that given by the
hyperbolic structure:
\begin{align*}
\begin{aligned}
T_0 = \rho_{\rm hyp}(T) &= \left(\begin{array}{rrrr}\thf&-\hf&\phm 1&\phm 0\\[1ex]\hf&\hf&1&0\\[1ex]1&-1&1&0\\[1ex]0&0&0&1\end{array}\right)
\qquad
U_0 = \rho_{\rm hyp}(U) &&= \left(\begin{array}{rrrr}\thf&-\hf&\phm \hf&\sqt\\[1ex]\hf&\hf&\hf&\sqt\\[1ex]\hf&-\hf&1&0\\[1ex]\sqt&-\sqt&0&1\end{array}\right)\\[2ex]
A_0 = \rho_{\rm hyp}(A) &= \left(\begin{array}{rrrr}1&0&0&\phm 0\\[1ex]0&-1&0&0\\[1ex]0&0&-1&0\\[1ex]0&0&0&1\end{array}\right)
\qquad
L_0 = \rho_{\rm hyp}(L) &&= \left(\begin{array}{rrrr}\phm 1&\phm 0&0&0\\[1ex]0&1&0&0\\[1ex]0&0&-\hf&\sqt\\[1ex]0&0&-\sqt&-\hf\end{array}\right)
\end{aligned}
\end{align*}

The first step of the process is to apply arbitrary small perturbations to the
generating matrices $T_0 \,,\, U_0 \,,\, A_0 \,,\, L_0$,
thus destroying the group relations.  We can however recover a representation
by instructing Newton's method to converge to matrices which {\it do} obey
the group relations.  Since the Newton iterative process is unlikely
to have awareness of the special nature of the original holonomy representation,
and since we already suspect from the Zariski tangent computation that the
variety has dimension $1$, we are not surprised when the process converges
to a representation $\rho'$ close to, but not conjugate to $\rho_{\rm hyp}$.

Let the images of $T \,,\, U \,,\, A \,,\, L$ under the new representation
$\rho'$ be $T_0' \,,\, U_0' \,,\, A_0' \,,\, L_0'$.  One readily observes
various features of the characteristic polynomials of these matrices.
In particular, the characteristic polynomials of $A_0' \,,\, L_0'$ are equal
to those of $A_0 \,,\, L_0$, respectively, and furthermore the characteristic
polynomials of $T_0' \,,\, U_0'$ are equal, of form
$$ 1  - (2+t)x + (2+2t)x^2 - (2+t)x^3 + x^4 \h.$$
It therefore seems reasonable to take \( t \) as a parameter for the
$1$--dimensional variety.  We choose a sequence
$t_1 \,,\, t_2 \,,\, \dots \,,\, t_n$ of rational values of $t$
(typically $n = 20$ suffices), and run the Newton process again $n$ times,
forcing it to converge to a sequence of representations
$\rho_1 \,,\, \rho_2 \,,\, \dots \,,\, \rho_n$, where $\rho_i(T)$
has characteristic polynomial
$1  - (2+t_i)x + (2+2t_i)x^2 - (2+t_i)x^3 + x^4$.

We aim to perform an interpolation of matrix entries of these $n$
representations, in order to obtain a {\it tautological representation}
encoding the entire variety, with generating matrices whose
entries are algebraic expressions in the parameter $t$.  However,
before this can be accomplished, it is necessary to conjugate each
representation $\rho_i$ into some ``normal'' form so that the
$\rho_i$ form an orderly sequence susceptible to interpolation.
This is technically the most demanding stage in the process;
it is described for representations of triangle groups in \cite{LT1}.

After achieving the normal forms for the $\rho_i$, we subject the
matrix entries to a lattice basis reduction algorithm, specifically
the LLL algorithm, to guess exact values for the $64$ matrix entries. 
The reader might justifiably feel some discomfort at the facts
(i) all computations to date have been numerical, and (ii) even after
the application of LLL we are still in the realm of guesswork.
However, as a final step we check that our conjectural tautological
representation is correct by formally verifying that it satisfies
the group relations, and that the ``base'' holonomy representation
$\rho_{\rm hyp}$ is given by a value of the parameter.

For an $r$--dimensional variety, the tautological representation arrived at
through this process is a representation of the group into a linear group
${\rm SL}(n,k)$, where the field $k$ is an algebraic extension of
$\Q(u_1, \dots u_r)$, the field of rational functions in indeterminates
$u_1, \dots u_r$ over $\Q$.  In the case of the Bianchi group ${\rm Bi}(3)$,
we have a single indeterminate \( u \) as the variety has dimension $1$;
moreover, we were able to conjugate so that $[k : \Q(u)] = 1$,
{\it i.e.} $k = \Q(u)$.
The resulting tautological representation
$\rho_u: {\rm Bi}(3) \to {\rm SL}(4,\Q(u))$ is as follows:
%
\begin{align*}
\rho_u(T) \;&=\; \left( \begin{array}{cccc}
1&0&0&0\\[2ex]
-\frac{2-u+u^2}{u^2} & \frac{1}{u} + u & -1 & \frac{1}{u} + u\\[2ex]
0&1&0&0\\[2ex]
0&0&0&1
\end{array}\right)\\[2ex]
\rho_u(U) \;&=\; \left(\begin{array}{cccc}
1+u&0&0&-u\\[2ex]
-\frac{2}{u+u^2} & \frac{1}{1+u} & \frac{1}{1+u} & -\frac{u}{1+u}\\[2ex]
-\frac{2+u+u^2}{u^2+u^3} & -\frac{1}{1+u} & \frac{1+u+u^2}{u+u^2} & \frac{1}{u+u^2}\\[2ex]
1&0&0&0
\end{array}\right)\\[2ex]
\rho_u(A) \;&=\; \left(\begin{array}{cccc}
1&0&0&0\\[2ex]
0&0&0&1\\[2ex]
-\frac{2-u+u^2}{u^2} & \frac{1}{u} + u & -1 & \frac{1}{u} + u\\[2ex]
0&1&0&0
\end{array}\right)\\[2ex]
\rho_u(L) \;&=\; \left(\begin{array}{cccc}
\frac{-1+u}{1+u} & -\frac{u^2}{1+u} & \frac{u}{1+u} & -\frac{u^2}{1+u}\\[2ex]
-\frac{2}{u+u^2} & \frac{1}{1+u} & \frac{1}{1+u} & -\frac{u}{1+u}\\[2ex]
-\frac{2+u+u^2}{u^2} & \frac{1}{u} & 0 & \frac{1}{u}\\[2ex]
-\frac{2}{u+u^2} & -\frac{u}{1+u} & \frac{1}{1+u} & \frac{1}{1+u}
\end{array}\right)
\end{align*}
We can think of $\rho_u$ as a family of representations parametrized by a complex
number $u \in \C \setminus \{ 0,-1 \}$.  The original holonomy representation
$\rho_{\rm hyp}$ is recovered up to conjugacy by setting $ u = 1 $.

\begin{prop} 
\begin{enumerate}
\item For any $u \in \C \setminus \{ 0,-1 \}$, $\rho_u: {\rm Bi}(3) \rightarrow {\rm SL}(4,\C)$ is a representation, i.e. the matrices $\rho_u(T), \rho_u(U), \rho_u(A)$ and $\rho_u(L)$ satisfy the relations in the Swan presentation of ${\rm Bi}(3)$.
\item When $u=1$, $\rho_1=\rho_{\rm hyp}$; for any $u, u' \in \C \setminus \{ 0,-1 \}$ with $u' \notin \{u, 1/u\}$, $\rho_u$ is not conjugate to $\rho_{u'}$.

\item For any $u \in {\rm U}(1) \setminus \{-1\}$, $\rho_u$ preserves a Hermitian form $H_u$, with signature $(3,1)$ when ${\rm Re} \, u > 1/4$ and $(4,0)$ when ${\rm Re} \, u < 1/4$.

\item For any $u \in {\rm U}(1) \setminus \{ \pm 1\}$, $\rho_u(T)$ and $\rho_u(U)$ are elliptic; in particular $\rho_u$ is not discrete or not faithful.

\item For any $u \in {\rm U}(1)$ with $1/4 < {\rm Re} \, u < 1$, $\rho_u ({\rm Bi}(3))$ is Zariski-dense in ${\rm SU}(3,1)$.
\end{enumerate}
\end{prop}

\Pf \begin{enumerate}
\item This follows from explicit computations best handled by formal computation software such as Mathematica or Maple. 
\item Likewise, using formal computation software one can readily check that when $|u|=1$, writing $s= {\rm Re} \, u$ and $t={\rm Im} \, u$, the four generators for $\rho_u ({\rm Bi}(3))$ preserve the Hermitian form given by the following matrix, where we denote $f(s,t)=-1-s+4s^2+i (-1+4s)t$:

$$
H_u = \left( \begin{array}{cccc}
-2+4s & f(s,t) & f(s,t) & f(s,t)  \\
\bar{f}(s,t)  & -2 + 4s & -1 & -1 \\
\bar{f}(s,t)  & -1 & -2 +4s & -1 \\
\bar{f}(s,t)  & -1 & -1 & -2+4s
\end{array}\right)
$$
The determinant of $H_u$ vanishes exactly when $s=1/4$; since the form has signature $(3,1)$ when $u=s=1$, $H_u$ has signature $(3,1)$ for all $u \in {\rm U}(1)$ with  ${\rm Re} \, u > 1/4$. Likewise, by testing a value in the other interval, we see that  the signature is $(4,0)$ when ${\rm Re} \, u < 1/4$.

\item By construction, $\rho_1=\rho_{\rm hyp}$. By a straightforward computation, the eigenvalues of $\rho_u(T)$ (as well as $\rho_u(U)$) are $(1,1,u,1/u)$, so the claim about non-conjugate $\rho_u,\rho_{u'}$ follows. 

\item We compute the dimension of the eigenspace of  $\rho_u(T)$ for the eigenvalue 1. Note that:

$$
\rho_u(T) - {\rm Id} =\; \left( \begin{array}{cccc}
0&0&0&0\\[2ex]
-\frac{2-u+u^2}{u^2} & \frac{1}{u} + u - 1 & -1 & \frac{1}{u} + u\\[2ex]
0&1&-1&0\\[2ex]
0&0&0&0
\end{array}\right)
$$
has rank 2, as its two non-zero lines are linearly independent. Therefore $\rho_u(T)$ is diagonalizable when $u \neq \pm 1$, so by the classification of isometries of complex hyperbolic space ${\rm H}_\C^3$ (see e.g. Theorem 3.4.1 of \cite{CG}), $\rho_u(T)$ is elliptic. The computation for $\rho_u(U)$ is similar. Recall that in the original representation, $\rho_{\rm hyp}(T)$ and  $\rho_{\rm hyp}(U)$ are parabolic, and in particular have infinite order. Now for $u \neq \pm 1$, either the elliptic isometry  $\rho_u(T)$ has finite order, in which case $\rho_u$ is not faithful, or it has infinite order, in which case $\rho_u$ has non-discrete image.

\item Let $G$ denote the Zariski-closure of $\rho_u ({\rm Bi}(3))$ in ${\rm SU}(3,1)$. In particular, $G$ is a closed subgroup of ${\rm SU}(3,1)$ with its classical topology. By Theorem 4.4.2 of \cite{CG}, either $G={\rm SU}(3,1)$ or $G$ has a global fixed point in $\overline{{\rm H}_\C^3}$ or preserves a totally geodesic submanifold of ${\rm H}_\C^3$. By the classification of totally geodesic subspaces (Proposition 2.5.1), in all these cases $G$ would be reducible, or contained in (a copy of) ${\rm SO}(3,1)$. It is easy to see that $\rho_u$ is irreducible for all $u$;  now since $\rho_u$ is not conjugate to $\rho_1$ for any $u \neq 1$, $\rho_u ({\rm Bi}(3))$ cannot be contained in a copy of ${\rm SO}(3,1)$ without violating rigidity of ${\rm Bi}(3)$ inside ${\rm SO}(3,1)$.
\end{enumerate}
\EPf


 

\section{Geometric structures and manifold subgroups}\label{structures}

\subsection{$\R$-Kleinian embeddings and tangent bundles}

The following is well known among experts, see for example section 2.5 of \cite{GKL} in dimension 2. It is interesting to note that the result holds in all dimensions, as opposed to the analogous interpretation of $\C$-Fuchsian representations which is only valid in dimension 2, by a coincidence of dimensions. We include a short proof, essentially the same as in \cite{GKL}, for the reader's convenience.

\begin{thm} Let $n \geqslant 2$, $M$ a complete hyperbolic $n$-manifold, $\rho_{hyp}: \Gamma = \pi_1(M) \longrightarrow {\rm SO}(n,1)$ the corresponding holonomy representation and $\iota$ the inclusion ${\rm SO}(n,1) \longrightarrow {\rm SU}(n,1)$. Then ${\rm H}_\C^n/\iota(\rho_{hyp}(\Gamma))$ is diffeomorphic to the tangent bundle $TM$. 
\end{thm}

\Pf $L={\rm H}^n_\R$ is a Lagrangian subspace of $X={\rm H}^n_\C$, so for any $p \in L$, $J \cdot T_p L$ is orthogonal to $T_p L$, where $J$ denotes the complex structure on $X$. In particular, $J \cdot T_p L \subset T_p X$ is the normal subspace to $L$ at $p$, hence $J$ induces an ${\rm SO}(n,1)$-equivariant diffeomorphism $TL \simeq NL$ (the latter denoting the normal bundle to $L$ in $X$). On the other hand, orthogonal projection onto $L$ induces a diffeomorphism $X \simeq NL$, also ${\rm SO}(n,1)$-equivariant. Therefore: $TM \simeq TL/\rho_{hyp}(\Gamma) \simeq NL/\iota(\rho_{hyp}(\Gamma)) \simeq X/\iota(\rho_{hyp}(\Gamma))$. \EPf

The same result holds, with the same proof, replacing $M$ by a hyperbolic $n$-orbifold and its tangent bundle by its tangent orbibundle. Roughly speaking, an orbifold is defined analogously to a manifold, but allowing certain points (the \emph{orbifold points}) to have neighborhoods diffeomorphic to a quotient of $\R^n$ by a finite isotropy group. Likewise, the tangent orbibundle of an orbifold is defined analogously to the tangent bundle of a manifold, but at the orbifold points the tangent space $R^n$ is replaced by its quotient under the isotropy group. See section 3.1 of \cite{C} for more details on tangent orbibundles.


\subsection{Deformations of manifold subgroups}\label{manifolds}

The small Bianchi groups are well known to contain knot and link groups as finite-index (torsion-free) subgroups. For example, Riley's original parametrization showed that the fundamental groups $\Gamma_8$ of the figure-eight knot complement is a subgroup of ${\rm Bi}(3)$ (with index 12), whereas ${\rm Bi}(1)$ contains the fundamental group $\Gamma_{\rm WL}$ of the Whitehead link complement (and ${\rm PGL}(2,\mathcal{O}_1)$ contains the Borromean rings group, see Section 9.2 of \cite{MR}).

The following result is Theorem 3.2 of \cite{BDL}; the dimension count appears in the proof of Lemma 3.6 therein:

\begin{thm}[\cite{BDL}] Let $M$ be an orientable complete finite volume hyperbolic manifold with fundamental group $\Gamma$ and $k$ cusps, and let $\rho_{\rm hyp}:\Gamma \longrightarrow {\rm SO}(3,1)$ be the holonomy representation of the complete hyperbolic structure. If $M$ is infinitesimally projectively rigid rel boundary, then $\rho_{\rm hyp}$ is a smooth point of ${\rm Hom}(\Gamma, {\rm SL}(4,\R))$, its conjugacy class  $[\rho_{\rm hyp}]$ is a smooth point of $\chi (\Gamma, {\rm SL}(4,\R))$, and ${\rm dim}_{[\rho_{\rm hyp}]} \, \chi (\Gamma, {\rm SL}(4,\R))=3k$.
\end{thm}

It was shown in \cite{HP} that $\Gamma_8$ and $\Gamma_{\rm WL}$ are infinitesimally projectively rigid rel boundary, giving the following:

\begin{cor} The character varieties $\chi (\Gamma_8,{\rm SL}(4,\R))$ and $\chi (\Gamma_{\rm WL},{\rm SL}(4,\R))$) are smooth at $[\rho_{{\rm hyp}}]$, with:
\begin{itemize}
\item ${\rm dim}_{[\rho_{{\rm hyp}}]} \, \chi (\Gamma_8,{\rm SL}(4,\R)) =3$ 
\item ${\rm dim}_{[\rho_{{\rm hyp}}]} \, \chi (\Gamma_{\rm WL},{\rm SL}(4,\R)) =6$
\end{itemize} 
\end{cor}

As before, by Theorem 2.2 of \cite{CLT1} (Theorem 1.4 above), one can replace ${\rm SL}(4,\R)$ with ${\rm SU}(3,1)$ in the above statements.



\raggedright

\begin{flushleft}
  \textsc{Julien Paupert\\
   School of Mathematical and Statistical Sciences, Arizona State University}\\
       \verb|paupert@asu.edu|
\end{flushleft}

\begin{flushleft}
  \textsc{Morwen Thistlethwaite\\
   Department of Mathematics, University of Tennessee Knoxville}\\
   \verb|mthistle@utk.edu|
\end{flushleft}

\end{document}